\documentclass{amsart}

\usepackage{amssymb}
\usepackage{amsmath}
\usepackage{comment}
\usepackage[all]{xy}
\usepackage{hyperref}

\newcommand{\R}{\mathbb R}
\newcommand{\N}{\mathbb N}
\newcommand{\Z}{\mathbb Z}
\newcommand{\ds}{\displaystyle}
\newcommand{\lie}{\mathcal{L}}
\renewcommand{\d}{\mathrm d}
\newcommand{\dbr}[1]{\left\{#1\right\}}
\newcommand{\ra}{\rightarrow}

\newcommand{\br}{\left[\, , \, \right]}
\newcommand{\brr}[1]{\left[#1\right]}
\newcommand{\X}{\ensuremath{\mathfrak{X}}}

\newcommand{\al}{\alpha}
\newcommand{\be}{\beta}

\newcommand{\m}{\mathrm{mod}\,}
\newcommand{\tr}{\mathrm{Tr}\,}

\newtheorem{theorem}{Theorem}           
\newtheorem{prop}[theorem]{Proposition}

\newtheorem{corol}[theorem]{Corollary}
\newtheorem{defi}[theorem]{Definition}

\newtheorem{obs}[theorem]{Remark}

\begin{document}

\title{Modular classes of Poisson-Nijenhuis Lie algebroids}
\author{Raquel Caseiro}
\address{CMUC, Department of Mathematics, University of Coimbra}
\email{raquel@mat.uc.pt} 

\begin{abstract}
The modular vector field of a Poisson-Nijenhuis Lie algebroid $A$ is
defined and we prove that, in case of non-degeneracy, this vector
field defines a hierarchy of bi-Hamiltonian $A$-vector fields. This
hierarchy covers an integrable hierarchy on the base manifold, which
may not have a Poisson-Nijenhuis structure.
\end{abstract}

\maketitle
\section{Introduction}

The relative modular class of a  Lie algebroid morphism was first
discussed by Gra\-bow\-ski, Marmo and Michor in \cite{GraMarMich}.
Kosmann-Schwarzbach and Weinstein in \cite{KosmannWeinstein} showed
that this relative class could be seen as a generalization of the notion
of modular class introduced by Weinstein in \cite{Weinstein}. In
\cite{LojaDamianou}, Damianou and Fernandes introduced the
modular vector field of a Poisson-Nijehuis manifold and showed that
it is intimately related with integrable hierarchies (see, also, the
alternative approach offered by Kosmann-Schwarzbach and Magri in
\cite{KosmannMagri0}). In this paper, we generalize this construction
and consider the modular vector field of a Poisson-Nijenhuis Lie algebroid.

Recall (see, e.g, \cite{KosmannMagri}) that a Nijenhuis operator
$N:A\ra A$ on a Lie algebroid $(A,\brr{~,~},\rho)$ allows us to
define a deformed Lie algebroid structure $A_N=(A,
\brr{~,~}_N,\rho\circ N)$ such that $N:A_N\to A$ is a Lie algebroid
morphism. Our first result states that the modular class of this
morphism has a canonical representative:

\begin{prop}
The relative modular class $N:A_N\to A$ is represented
by $\d_A \tr N$.
\end{prop}

Let us assume now that $A$ is equipped with a Poisson structure
$\pi$ compatible with $N$. Then we can define two Lie algebroid
structures on $A^*$, namely $(A^*,\brr{~,~}_\pi,
\rho\circ\pi^\sharp)$ obtained by dualization from $\pi$, and
$A_{N^*}^*=(A^*, [~ , ~]_{N\pi}, \rho\circ N\pi^\sharp)$ obtained
from the first one by deformation along $N^*$. Again,
$N^*:A_{N^*}^*\to A^*$ is a Lie algebroid morphism and, by the
proposition above, its relative modular class has the canonical
representative $X_{(N,\pi)}:=\d_\pi(\tr N)$, which we will call the
\emph{modular vector field} of the Poisson-Nijenhuis Lie algebroid
$(A,\pi,N)$. When $A=TM$ we recover the construction of Damianou and
Fernandes \cite{LojaDamianou} up to a factor of $1/2$ (for the same
reason that the modular class associated with a Poisson manifold
differs from the modular class of its cotangent Lie algebroid by a
factor of $1/2$).

The modular vector field $X_{(N,\pi)}$ of the Poisson-Nijenhuis Lie
algebroid $(A,\pi,N)$ is always a $\d_{N\pi}$-cocycle. If $N$ is
non-degenerated it is also a $\d_{N\pi}$-coboundary. In this case we have
the following generalization of a result of Damianou and Fernandes
for a Poisson-Nijenhuis manifold:

\begin{theorem}
\label{thm:main:intr}
Let $(A,\pi,N)$ be a Poisson-Nijenhuis Lie algebroid with $N$
a non-degenerated Nijenhuis operator. Then there exists a
hierarchy of $A$-vector fields
\[ X_{(N,\pi)}^{i+j}=N^{i+j-1}X_{(N,\pi)}=\d_{N^i\pi}h_j=\d_{N^j\pi}h_i,
\quad (i,j\in\Z) \]
where
\[ h_0=\ln (\det N) \mbox{ and } \ds h_i=\frac{1}{i}\tr{N^i},\quad (i\neq 0).\]
\end{theorem}

The hierarchy of flows on $A$, given by this theorem, covers a hierarchy
of (ordinary) multi-Hamiltonian flows on the base manifold $M$. Although
the hierarchy on $A$ is generated by a Nijenhuis operator, it may happen
that the base hierarchy is not generated by one. We will see that this is
precisely the case for the $A_n$-Toda lattice. This gives a new explanation
for the existence of a hierarchy of Poisson structures and flows associated
with a bi-Hamiltonian system which may not have a Nijenhuis operator.

This paper is organized as follows. In Section 2, we present the necessary
background on Poisson-Nijenhuis Lie algebroids. In Section 3, we introduce
the modular vector field of a Poisson-Nijenhuis Lie algebroid, state
its basic properties, and prove Theorem \ref{thm:main:intr}. Section 4 is
concerned with integrable hierarchies and discusses the example of the $A_n$-Toda
lattice.

\section{Poisson-Nijenhuis Lie Algebroids}
In this section we will recall some basic facts about Nijenhuis
operators and Poisson structures on Lie algebroids which we will
need later. A general reference for Lie algebroids is the book by
Cannas da Silva and Weinstein \cite{CaWe}. Nijenhuis operators are
discussed in detail in the article by Kosmann-Schwarzbach and Magri
\cite{KosmannMagri}, while PN-structures on Lie algebroids are
discussed by Kosmann-Schwarzbach in \cite{Kosmann} and by Grabowski
and Urbanski in \cite{GraUrb}.

\subsection{Cartan calculus on Lie algebroids}
Let us recall that a Lie algebroid is a kind of generalized tangent
bundle, which carries a generalized Cartan calculus.

First, for any Lie algebroid $(A,\br,\rho)$ we have a complex
of \emph{$A$-differential forms} $\Omega^k(A):=\Gamma(\wedge^kA^*)$
with the differential given by:
\begin{align*} \d_A \omega(X_0,\dots,X_k):=\sum_{i=0}^n(-1)^i\rho(X_i)\cdot \omega(X_0,\ldots, \hat{X_i},\ldots, X_k)&\\
+\sum_{0\leq i\leq j\leq k} (-1)^{i+j}\omega \left([X_i,X_j]_A,X_0,\ldots, \hat{X_i},\ldots, X_k \right)&,
\end{align*}
where $X_0,\ldots, X_k\in\Gamma(A)$.
The corresponding \emph{Lie algebroid cohomology} is denoted by $H^\bullet(A)$.

Dually, the space of \emph{$A$-multivector fields}
$\X^\bullet(A)=\bigoplus_{k\in\Z}\X^k(A):=\bigoplus_{k\in\Z}\Gamma(\wedge^k
A)$ carries a super-Lie bracket $[~,~]_A$, extending the Lie bracket
on $\Gamma(A)$, and satisfying the following super-commutation,
super-derivation and super-Jacobi identities:
\begin{align*}
     &[P,Q]=-(-1)^{(p-1)(q-1)} [P,Q]\\
     &[P,Q\wedge R]=[P,Q]\wedge R+(-1)^{(p-1)q}Q\wedge[P,R]\\
     &(-1)^{(p-1)(r-1)} [P,[Q,R]]+ (-1)^{(q-1)(p-1)} [Q,[R,P]]
      +(-1)^{(r-1)(q-1)} [R,[P,Q]]=0
\end{align*}
where $P\in\X^p(A)$, $Q\in\X^q(A)$ and $R\in\X^r(A)$. The triple $(\X^\bullet(A),[~,~]_A,\wedge)$
is a \emph{Gerstenhaber algebra}.

If $X\in\Gamma(A)$ and $\omega\in\Omega^k(A)$ the \emph{Lie derivative} of $\omega$ along
$X$ is the $A$-differential form $\lie_X\omega\in\Omega^k(A)$ defined by
\[ \lie_X\omega:=\d_A i_X\omega+i_X\d_A\omega,\]
where $i_X:\Omega^k(A)\ra\Omega^{k-1}(A)$ is   defined by
\begin{align*}
i_X\eta(X_1,\ldots X_{k-1})&:=\eta(X,X_1,\ldots,X_k), \quad X_1,\ldots, X_{k-1}\in\Gamma(A),
\end{align*}
for $k> 1$. If $k=1$, then $i_X\eta:=\eta(X)$ and, for $k\leq 0$ we say that  $i_X\eta=0$.

A morphism $\phi:A\to B$ (over the identity) of Lie algebroids over $M$ induces
by transposition a chain map of the complexes of differential forms:
$$\phi^\ast:(\Omega^k(B),d_B)\ra(\Omega^k(A),\d_{A}).$$
Hence, we also have a well defined map at the level of cohomology, which
we will denote by the same letter $\phi^\ast:H^\bullet(B)\ra H^\bullet(A)$.

\subsection{Nijenhuis operators}

Let $(A,\br,\rho)$ be a Lie algebroid over a manifold $M$. Recall that
a \emph{Nijenhuis operator} is a bundle map $N:A\ra A$ (over the identity) such
that the induced map on the sections (denoted by the same symbol $N$) has
vanishing torsion:
\begin{equation}
\label{eq:Nijenhuis} T_N(X,Y):=N[X,Y]_N-[NX,NY]=0, \quad  X,Y\in
\Gamma(A),
\end{equation}
where $\br_N$ is defined by
$$
[X,Y]_N:=[NX,Y]+[X,NY]-N[X,Y],\quad X,Y\in \Gamma(A).
$$
Let us set $\rho_N:=\rho\circ N$. For a Nijenhuis operator $N$,
one checks easily that the triple $A_N=(A,\br_N,\rho_N)$ is a new
Lie algebroid, and then $N:A_N\ra A$ is a Lie algebroid morphism.

Since $N$ is a Lie algebroid morphism, its transpose gives a chain
map of the com\-ple\-xes of dif\-fe\-ren\-tial forms
$N^\ast:(\Omega^k(A),d_A)\ra(\Omega^k(A_N),\d_{A_N})$. Hence we also
have a map at the level of algebroid cohomology
$N^\ast:H^\bullet(A)\ra H^\bullet(A_N)$.

\subsection{Poisson structures on Lie algebroids}

Let  $\pi\in\X^2(A)$ be a bivector on the Lie algebroid
$(A,\br,\rho)$ and denote by $\pi^\sharp$ the usual  bundle map
\[
\begin{array}{llll}
  \pi^\sharp: &  A^\ast & \longrightarrow &A\\
  & \al &\longmapsto  & \pi^\sharp(\al)=i_\al \pi. \\
\end{array}
\]
We say that $\pi$ defines a \emph{Poisson structure on $A$}
if $[\pi,\pi]_A=0$. In this case, the  bracket on the sections of $A^\ast$ defined by
\begin{align*}
[\al,\be]_{\pi}&=\mathcal{L}_{\pi^\sharp\al}\be-\mathcal{L}_{\pi^\sharp\be}\al-\d_A\left(\pi(\al,\be)\right),
\quad \al,\be\in\Gamma(A^\ast),
\end{align*}
is a Lie bracket and $(A^\ast, \br_{A^\ast}, \rho\circ\pi^\sharp)$ is a Lie algebroid.
The differential of this Lie algebroid is given by
$\d_\pi X=[\pi , X]_A$, $X\in\Omega(A^\ast)$, and the pair
$(A,A^\ast)$ is a special kind of Lie bialgebroid, called a \emph{triangular Lie bialgebroid}.

\subsection{Poisson-Nijenhuis Lie algebroids} The basic notion to be used
in this paper is the following:

\begin{defi}
A \textbf{Poisson-Nijenhuis Lie algebroid} (in short, a PN-algebroid)
is a Lie algebroid $(A,\br_A,\rho)$ equipped with a Poisson structure $\pi$ and
a Nijenhuis operator $N$ which are \emph{compatible}.
\end{defi}

The compatibility condition between $N$ and $\pi$ means that:
\[ \br_{N\pi}=\br_{N^\ast},\]
where $\br_{N\pi}$ is the Lie bracket defined by the bivector field
$N\pi\in\X^2(A)$, and $\br_{N^\ast}$ is the Lie bracket obtained from
the Lie bracket $\br_\pi$ by deformation along the Nijenhuis tensor
$N^\ast$.

As a consequence, $N\pi$ defines a new Poisson structure on $A$, compatible with
$\pi$:
\[ [\pi,N\pi]_A=[N\pi,N\pi]_A=0,\]
and one has a commutative diagram of morphisms of Lie algebroids:
\[
\xymatrix{
(A^\ast,[\cdot,\cdot]_{N\pi})\ar[rr]^{{N}^*}\ar[dd]_{\pi^\sharp}\ar[ddrr]^{N\pi^\sharp}
&&
(A^\ast,[\cdot,\cdot]_{\pi})\ar[dd]^{\pi^\sharp}\\
\\
(A,[\cdot,\cdot]_{N})\ar[rr]^{{N}}&&(A,[\cdot,\cdot]_A) }
\]

In fact, we have a whole hierarchy ${N^k\pi}$ ($k\in\N$) of pairwise compatible Poisson
structures on $A$.

\section{Modular class of a Poisson-Nijenhuis Lie algebroid}

In this section we will state and prove our main results.

\subsection{Modular class of a Lie algebroid}
Let $(A, \br, \rho)$ be a Lie algebroid over the manifold $M$. For
simplicity we will assume that both $M$ and $A$ are orientable, so
that there exist non-vanishing sections
$\eta\in\X^{\mathrm{top}}(A)$ and $\mu\in\Omega^{\mathrm{top}}(M)$.

The \emph{modular form} of the Lie algebroid $A$ with respect to $\eta\otimes\mu$
is the $A$-form $\xi_A\in\Omega^1(A)$, defined by
\begin{equation}
    \langle \xi_A,X\rangle \eta\otimes \mu=\lie_X\eta\otimes
    \mu+\eta \otimes\lie_{\rho(X)}\mu, \quad X\in\Gamma(A).
\end{equation}
This is a 1-cocycle of the Lie algebroid cohomology of $A$. If one makes
a different choice of sections $\eta'$ and $\mu'$, then
$\eta'\otimes\mu'=f\eta\otimes \mu$, for some non-vanishing smooth function
$f\in C^\infty(M)$. One checks easily that the modular form $\xi'_A$ associated
with this new choice is given by:
\begin{equation}
\label{eq:change:section}
\xi'_A=\xi_A-\d_A\log|f|,
\end{equation}
so that the cohomology class $[\xi_A]\in H^1(A)$ is independent of the choice of
$\eta$ and $\mu$. This cohomology class is called the \emph{modular class} of $A$
and we will denoted it by $\m A:=[\xi_A]$.

\begin{prop}
\label{modN} Let $N$ be a Nijenhuis operator on a Lie algebroid $A$
and fix non-vanishing sections $\eta$ and $\mu$ as above. The
modular form $\xi_{A_N}$ of the Lie algebroid $A_N$ and the modular
form $\xi_A$ of $A$ are related by:
\begin{equation}
 \xi_{A_N}=\d_A(\tr N)+N^\ast\xi_A.
\end{equation}
\end{prop}

\begin{proof}
Around any point, we can always choose a local base $\{e_1,\ldots,
e_r\}$ of sections of $A$ and local coordinates $(x_1,\ldots, x_n)$
of $M$ such that $\eta=e_1\wedge\ldots\wedge e_r$ and
$\mu=dx_1\wedge\ldots\wedge dx_n$. In these coordinates, we have the
following expressions for the anchor $\rho$ and the Nijenhuis
operator $N$:
\[
\rho(e_i)=\sum_{u=1}^n p_i^u\,\frac{\partial}{\partial x_u} \text{ and }
N(e_i)=\sum_{j=1}^rN_i^je_j,\quad (i=1,\dots,r).
\]
Now, for $i=1,\dots,r$, we compute:
$$
\lie_{e_i}\eta=\lie_{e_i}(e_1\wedge\ldots\wedge e_r)=\sum_{k,j=1}^r
\left[N_i^kC_{kj}^j-\sum_{u=1}^n\left(\rho_j^u\frac{\partial
N_i^j}{\partial x_u}+ \rho_i^u\frac{\partial N_j^j}{\partial
x_u}\right)\right]\eta
$$
and
$$
\lie_{\rho_N(e_i)}\mu=\lie_{\rho\circ
N(e_i)}(dx_1\wedge\ldots\wedge dx_n)=\sum_{k=1}^r\sum_{u=1}^n
\left(\frac{\partial N_i^k}{\partial x_u}\rho_k^u+\frac{\partial
\rho_k^u}{\partial x_u} N_i^k\right)\mu.
$$
So
\begin{eqnarray*}
\xi_A^N(e_i)&=&\sum_j\left(\sum_k N_i^jC_{jk}^k+\sum_u\left(\frac{\partial \rho_j^u}{\partial x_u}N_i^j
+\frac{\partial N_i^j}{\partial x_u}\rho_i^u\right)\right)\\
&=&
\sum_j\left(\sum_k N_i^jC_{jk}^k+\sum_u\frac{\partial \rho_j^u}{\partial x_u}N_i^j\right)+\d_A\tr N(e_i)\\
&=& \left( N^\ast\xi_A+\d_A\tr N\right)(e_i).
\end{eqnarray*}

By linearity this holds for any section of $A$, so the result
follows.
\end{proof}

The theorem shows that the $A_N$-form $\xi_{A_N}-N^\ast\xi_A=\d_A\tr N$ is independent of the
choice of section of $\eta\otimes\mu\in\X^{\mathrm{top}}(A)\otimes\Omega^{\mathrm{top}}(M)$.

\begin{obs}
This can also be checked directly using relation (\ref{eq:change:section}): If
$\xi_A$ and $\xi_{A_N}$ are the modular forms associated to the choice $\eta\otimes\mu$,
$\xi'_A$ and $\xi'_{A_N}$ are the modular forms associated with another choice
$f\eta\otimes\mu$, then:
\begin{align*}
 \xi'_{A}&=\xi_{A}-\d_{A}\ln |f|,\\
 \xi'_{A_N}&=\xi_{A_N}-\d_{A_N}\ln |f|.
\end{align*}
For any function $g\in C^{\infty}(M)$, we have $\d_{A_N}g=N^\ast d_A g$, so
it follows from these relations that:
\[
\xi'_{A_N}-N^\ast\xi'_A=\xi_{A_N}-N^\ast\xi_A.
\]
\end{obs}

Recall (see \cite{GraMarMich,KosmannWeinstein}) that for any
Lie algebroid morphism over the identity
$\phi:(A,\br,\rho_A)\ra (B,\br_B, \rho_B)$ one defines its
\emph{relative modular class} to be the cohomology class
$\mathrm{mod}^\phi(A,B)\in H^1(A)$ given by:
\begin{equation}
    \mathrm{mod}^\phi(A,B):=\m A-\phi^\ast \m B.
\end{equation}
Therefore we have the following immediate corollary of Proposition \ref{modN}:

\begin{corol}
\label{coro:rel}
The relative modular class of the algebroid morphism $N:A_N\ra A$ is
a $A_N$-cohomology class with canonical representative the $A_N$-form
$\d_A \tr N$.
\end{corol}

Note that the class $[\d_A \tr N]\in H^1(A_N)$ maybe non-trivial: in general,
the differentials $\d_A$ and $\d_{A_N}$ will be distinct.

\subsection{Modular class of a Poisson-Nijenhuis Lie algebroid}

Now we consider a Poisson-Nijenhuis Lie algebroid
$\left(A,\pi,N\right)$. Then $N^\ast$ is a Nijenhuis operator of the
dual Lie algebroid $(A^\ast, \brr{~,~}_\pi, \rho\circ\pi^\sharp)$
and, by Corollary \ref{coro:rel}, its relative modular class has the
canonical representative $\d_{\pi}(\tr N^\ast)$, so that:
$$
\mathrm{mod}^{N^\ast}(A^\ast_{N^\ast},A^\ast)=[\d_{\pi}(\tr N^\ast)]=[\d_{\pi}(\tr N)].
$$

\begin{defi}
  The \textbf{modular vector field} of the Poisson-Nijenhuis Lie algebroid
  $(A,\pi, N)$ is defined by
    $$X_{(N,\pi)}=\xi_{A^\ast_{N^\ast}}-N\xi_{A^\ast}=\d_{\pi}(\tr N)\in\X(A).$$
\end{defi}

Notice that the modular vector field of a PN-algebroid is a
$\d_{N\pi}$-cocycle and we have a generalization to PN-algebroids of
the results obtained in \cite{LojaDamianou, KosmannMagri0} for a
Poisson manifold.

\begin{prop}\label{modhier}
Let $\left(A,\pi,N\right)$ be a PN-algebroid. Then
\begin{equation*}
N^kX_{(N,\pi)}
=\frac{1}{k-i+1}\,X_{(N^{k-i+1},\,N^i\pi)}, \quad i<k\in\N.
\end{equation*}
\end{prop}

\begin{proof}
The operator $N$ is Nijenhuis  so it satisfies the identity

\begin{equation}\label{6}
   k N^{\ast\, k}\d_A\tr N=\d_A \tr N^k, \quad k\in\N.
\end{equation}

 Now simply observe that
\begin{align*}
N^{k}X_{(N,\pi)}&= N^{k}\d_{\pi}(\tr N)=N^{k}\, [\pi,\tr N]\\
                &=-N^i\,\pi^\sharp(N^{\ast\,k-i}\, \d_A \tr N)
                =-\frac{1}{k-i+1}\,(N^i\pi)^\sharp\left(\,\d_A\tr N^{k-i+1}\right)\\
                  &=\frac{1}{k-i+1}\,\d_{N^i\pi}\left(\tr N^{k-i+1}\right)
                 =\frac{1}{k-i+1}\,X_{(N^{k-i+1},\,N^i\pi)}.
\end{align*}
\end{proof}

As an immediately consequence  we have:

\begin{corol} The modular vector field determines a hierarchy of vector
fields
\begin{equation}\label{eq:A:hierarchy:one}
X^{k}_{(N,\pi)}=N^{k-1}X_{(N,\pi)}=\,\d_{N^{k-i}\pi}h_{i}=\,\d_{N^{i-1}\pi}h_{k-i+1},
\quad (i\leq k\in\N, k>1),
\end{equation}
where $h_i=\frac{1}{i}\tr N^i$, $i\in\N$.
\end{corol}

\begin{obs}
The basic formula (\ref{6}) is well-known to people working in
integrable systems (see for instance \cite{Magri}). Also, special
versions of Proposition \ref{modhier} and its corollary can be found
in \cite{LojaDamianou, KosmannMagri0}, which were sources of
inspiration for this work.
\end{obs}

If $N$ is non-degenerated then the modular vector field  also
defines a negative hierarchy.

\begin{theorem}
\label{thm:main} Let $(A,\pi,N)$ be a Poisson-Nijenhuis Lie
algebroid with $N$ a non-degenerated Nijenhuis operator. Then the
modular vector field $X_{(N,\pi)}$ is a $\d_{N\pi}$-coboundary  and
determines a hierarchy of vector fields
\begin{equation}
\label{eq:A:hierarchy}
X_{(N,\pi)}^{i+j}=N^{i+j-1}X_{(N,\pi)}=\d_{N^i\pi}h_j=\d_{N^j\pi}h_i,
\quad (i,j\in\Z)
\end{equation}
where
\begin{equation}
\label{eq:Hamiltonians} h_0=\ln (\det N) \text{ and } \ds
h_i=\frac{1}{i}\tr{N^i},\quad (i\neq 0).
\end{equation}
\end{theorem}

\begin{proof}
For any integer $k$, $N^k$ is a non-degenerated Nijenhuis operator
and satisfies the identity
\begin{equation}\label{nijdet}
k N^{\ast \,k} \d_A(\ln\det N)=\d_A(\tr N^k).
\end{equation}
It follows that
$$
{\ds X_{(N,\pi)}=-\pi^\sharp\left( \d_A \tr N \right) =-\pi^\sharp
\,N^\ast\d_A\left(\ln\det N\right) =\d_{N\pi}(\ln\det N)}.
$$

We also have
\begin{align*}
N^{-1}X_{(N,\pi)}&=X_{(N,N^{-1}\pi)}=-N^{-1}\pi^\sharp(\d_A \tr N)
                                    =-\pi^\sharp(\d_A\ln\det N)\\
                                    &=- \pi^\sharp\left( N^\ast \d_A \tr N^{-1} \right)
                                    =\d_{N\pi}\tr N^{-1}=X_{(N^{-1},N\pi)}.
\end{align*}
The expressions for the hierarchy now follow from identity
(\ref{nijdet}) and Proposition \ref{modhier} applied to $N$ and to
$N^{-1}$.
\end{proof}

\begin{obs}
In case of degeneracy of $N$ one can always consider a
non-degenerated operator of the form $N+\lambda I$, $\lambda$
constant. Although this Nijenhuis operator does not define exactly
the same hierarchy as $N$,  the  traces of its powers generate the
same  algebra as the functions $h_i$. I thank an anonymous referee
for this remark.
\end{obs}

\section{Integrable hierarchies and PN-Algebroids}

Let $(A, \br,\rho)$  be a Lie algebroid over a manifold $M$ and
consider the fiberwise linear Poisson structure $\{~,~\}_A$ on the
dual bundle $A^*$. A vector field $X\in\X(A)$ defines, by
evaluation, the function $f_X: A^*\to \R$, which is linear on the
fibers:
\begin{equation*}
  f_X(\alpha(x))=\langle \alpha(x) , X(x)\rangle, \quad \al\in\Gamma(A^*).
\end{equation*}
We denote by $X_{f_X}$ its the Hamiltonian vector field. It is easy
to check (see~\cite{Loja, GraUrb, GraUrb2}) that:
\begin{enumerate}
\item[(a)] The assignment $X\mapsto f_X$ defines a Lie algebra
homomorphism
\[ (\X(A),[~,~])\to (C^\infty(A^*),\{~,~\}_A);\]
\item[(b)] Denoting by $q: A^*\to M$ the projection,
$X_{f_X}$ is $q$-related to $\rho(X)$:
\[ q_*X_{f_X}=\rho(X).\]
\end{enumerate}

So the vector field $X_{f_X}$ is  the  lift of the vector field
$\rho(X)$ to $A^*$ (associated with the  derivation on $A^\ast$,
$D_X^\ast\alpha=\mathcal{L}_X\alpha$).

One can also define $X^A$, the lift of $\rho(X)$ to $A$ (associated
with the derivation on $A$, $D_X =\brr{X,-}$) by
 \begin{align*}
 &X^A(f_\alpha)=f_{\mathcal{L}_X\alpha},\quad\alpha\in\Gamma(A^*)\\
&X^A(g\circ p)=X(g)\circ p,\quad g\in C^\infty(M),
\end{align*}
where $f_\alpha:A\ra\R$ is linear function defined by evaluation by
$\alpha\in \Gamma(A^\ast)$  and
 $p:A\ra M$ is the
projection of the vector bundle.

\subsection{Integrable Hierarchies on Lie algebroids}
Let  $\pi\in\X^2(A)$ be a Poisson structure on a Lie algebroid $A$
over a manifold $M$. This bivector field determines a bundle map
$\pi^\sharp:A^*\to A$, and for the Lie algebroid structure on $A^*$,
we have that $\pi^\sharp$ is a Lie algebroid morphism.

\begin{defi}
Let $f\in C^\infty(M)$. Its \textbf{Hamiltonian vector field}
$X_f\in\X(A)$ is the vector field:
\[ X_f:=\pi^\sharp \d_A f.\]
\end{defi}

The Poisson structure $\pi$ on $A$ covers a (ordinary) Poisson
structure $\pi_M$ on the base manifold $M$  which is defined by
${\ds \pi_M^\sharp=\rho\circ\pi^\sharp\circ\rho^\ast}$, i.e.
\begin{equation*}
\dbr{f,g}_{\pi_M}=\langle \d_A f,\d_{\pi}g\rangle=\pi(\d_A f,\d_A
g),\quad f,g\in C^\infty(M).
\end{equation*}

Moreover, for any $f\in C^\infty(M)$ the Hamiltonian vector field
$X_f$ on $A$ covers the (ordinary) Hamiltonian vector field on $M$
associated with $f$:
$$\rho(X_f)=\rho\circ\pi^\sharp(\d_A f)=
\rho\circ\pi^\sharp\circ\rho^*(\d f)=\pi_M^\sharp(\d f).$$

Considering $\{~, ~\}_{A^*}$, the linear Poisson bracket on $A$
defined by the Lie algebroid ${\ds (A^*, \brr{~,~}_\pi,
\rho\circ\pi^\sharp)}$, notice that the Hamiltonian vector field
$X_{f_{\d_A f}}$ is the  lift of $\rho(X_f)$ to $A$, associated with
the derivation ${\ds D_{X_f}=[X_f, ~ ]}$ and $X_{f_X}$ is the
 lift of the vector field $\rho(X)$ to $A^*$, associated
with the derivation $D_{\d _A f}=D_{X_f}^*$.

\begin{defi}
Given $\pi$ a Poisson structure on $A$, we say that two functions
$g,f\in C^\infty(A)$ \textbf{$\pi$-commute} if $\dbr{f,g}_{A^*}=0$.
A first integral of  a vector field $X$ is  a first integral of
$X^A$, the
 lift of $\rho(X)$ to $A$, and we say that $X$ is  \textbf{integrable} if $X^A$
Liouville integrable.
\end{defi}

Two  first integrals of the vector field $\rho(X)$ may not
$\pi_M$-commute but their pull-backs always $\pi$-commute, because
basic functions always commute with respect to the linear Poisson
bracket defined on the dual of a Lie algebroid. In particular,
we have:

\begin{prop}
Let ${\ds f\in C^\infty(M)}$. The first integrals of the Hamiltonian vector
field
\[ X_{f}=\pi^\sharp \d_A f,\]
are the functions which $\pi$-commute with ${\ds f_{\d_A f}}$. In particular,
evaluations of sections of $A^*$ which commute with ${\ds \d_A f}$ and pull-backs of
functions which $\pi_M$-commute with $f$ are first integrals of ${\ds X_f}$.
\end{prop}

Let $N$ be a Nijenhuis operator on $A$ compatible with $\pi$. The sequence
of Poisson structures $\pi_k=N^k\pi$ covers a sequence of Poisson structure on $M$:
\begin{equation*}
    {\pi^{k}_M}^\sharp=\rho\circ N^k\pi\circ\rho^*,
\end{equation*}
but, as an example below shows, these Poisson tensors may not be related
by a Nijenhuis operator on $M$. We also have a sequence of linear Poisson brackets on $A$:
\begin{equation*}
    \dbr{ {~} , ~ }_{A^*}^k
\end{equation*}
and, given a section $\al$ on $A^*$, a sequence of Hamiltonian vector fields
\begin{equation*}
    X^k_{f_\alpha}=\{f_\alpha, - \}_{A^\ast}^k.
\end{equation*}

A bi-Hamiltonian vector field
\[
X=\pi_0^\sharp(\d_A h_1)=\pi_1^\sharp(\d_A h_0)
\]
defines a hierarchy of multi-Hamiltonian vector fields on $A$:
\[
X_{k+i}=\pi_k^\sharp(\d_A h_i)=\pi_i^\sharp(\d_A h_k).
\]
This hierarchy, on one hand, covers a hierarchy of Hamiltonian
vector fields on $M$
\begin{equation*}
    \rho(X_{i+k})={\pi^{k}_M}^\sharp\d h_i={\pi^{i}_M}^\sharp\d
    h_k
\end{equation*}
and, on the other hand, is associated with the hierarchy of
 lifts
\[
\rho(X_{i+k})^A=X^0_{f_{\d_A h_{i+k}}}=
X^i_{f_{\d_A h_{k}}}=X^k_{f_{\d_A h_{i}}}.
\]

\subsection{Covering Integrable Hierarchies}

We can try to apply our main result (Theorem \ref{thm:main}) to
obtain an integrable hierarchy on a Lie algebroid. However, one
observes that in Theorem \ref{thm:main} the Hamiltonian functions,
which are first integrals of the vector fields in the hierarchy, are
all basic functions, i.e., are pull-backs of functions on the base.
Hence, in general, they will not provide a complete set of first
integrals. However, there is yet another connection with (classical)
integrable systems, due to the following theorem:

\begin{theorem}
\label{thm:covering:hier} Let $(A,\pi,N)$ be a Poisson-Nijenhuis Lie
algebroid with $N$ a non-degenerated Nijenhuis operator. Then the
the modular vector field $X_{(N,\pi)}$ covers a bi-Hamiltonian
vector field on $M$, and the associated hierarchy
(\ref{eq:A:hierarchy}) of $A$-vector fields covers a (classical)
hierarchy of vector fields on $M$. This hierarchy is given by:
\begin{equation}
\label{eq:hierarchy} X_{i+j}=-\pi_i^\sharp\d h_j=-\pi_j^\sharp\d
h_i\quad (i,j \in \Z)
\end{equation}
where $\pi_j$ are Poisson structures on $M$ and $h_i$ are the functions
given by (\ref{eq:Hamiltonians}).
\end{theorem}

Although we have a hierarchy of  modular vector fields
$X_{(N^k,\pi)}$ generated by the Nijenhuis operator $N$, generally
the covered hierarchy of bi-Hamiltonian vector fields on $M$ is not
 generated by any Nijenhuis operator. This is illustrated by the
 next example.


\subsection{The classical Toda lattice}

The classical Toda lattice was already considered in \cite{Agrotis},
using specific properties of this system. We use our general
approach to show how one can recover those results and explain some
of those formulas.

\subsubsection{Toda lattice in physical coordinates}

The Hamiltonian defining the Toda lattice is given in canonical
coordinates $(p_i,q_i)$ of $\R^{2n}$ by
\begin{equation}
\label{eq:Toda:An} h_2(q_1,\dots,q_n,p_1,\dots,p_n)=
\sum_{i=1}^n\,\frac{1}{2}\,p_i^2+\sum _{i=1}^{n-1}\,e^{q_i-q_{i+1}}.
\end{equation}
For the integrability of the system we refer to the classical paper
of Flaschka \cite{Flaschka1}.

Let us recall the bi-Hamiltonian structure given in
\cite{Fernandes1}. The first Poisson tensor in the hierarchy is the
standard canonical symplectic tensor, which we denote by
$\widetilde{\pi}_0$, so that
\[ \{q_i,p_j\}_0=\delta_{ij},\]
while the second Poisson tensor $\widetilde{\pi}_1$ is determined by
the relations
\begin{align*}
\{q_i,q_j\}_1&=-1,\qquad (i<j)\\
\{q_i,p_j\}_1&=p_i\delta_{ij},\\
\{p_j,p_i\}_1&=e^{q_i-q_{i+1}}\delta_{j,i+1}.
\end{align*}
Then setting $h_1=p_1+p_2+\cdots+p_n$, we obtain the bi-Hamiltonian
formulation:
\[ \widetilde{\pi}_0^\sharp\d h_2=\widetilde{\pi}_1^\sharp \d h_1.\]
If we set, as usual,
\[ N:=\widetilde{\pi}_1^\sharp\circ(\widetilde{\pi}_0^\sharp)^{-1},\]
then a small computation  gives the following multi-Hamiltonian
formulation:

\begin{prop}\label{propTodaMulti}
The Toda hierarchy admits the multi-Hamiltonian formulation:
\[ \widetilde{\pi}_j^\sharp \d h_2=\widetilde{\pi}_{j+2}^\sharp \d h_0,\]
where $h_0=\frac{1}{2}\log({\det N})$ and $h_2$ is the original
Hamiltonian (\ref{eq:Toda:An}).
\end{prop}

\subsubsection{Toda lattice in Flaschka coordinates}

Let us recall the Flaschka coordinates
$(a_1,\dots,a_{n-1},b_1,\dots,b_n)$ where:
\begin{align*}
b_i&=p_i, \quad (i=1,\dots,n)\\
a_i&=e^{q_i-q_{i+1}}, \quad (i=1,\dots,n-1)
\end{align*}
In these new coordinates there is no recursion operator anymore
(this is a singular change of coordinates, where we loose one degree
of freedom). Nevertheless, the multi-Hamiltonian structure does
reduce (\cite{Fernandes1}). One can then compute the modular vector
fields of the reduced Poisson tensors $\overline{\pi}_j$ relative to
the standard volume form:
\[
\mu=\d a_1\wedge\cdots \wedge\d a_{n-1}\wedge\d
b_1\wedge\cdots\wedge\d b_n.
\]

It turns out that the modular vector fields $X^j_\mu$ are
Hamiltonian vector fields with Hamiltonian function
\[ h=\log(a_1\cdots a_{n-1})+ j \log(\det(L)), \]
where $L$ is the Lax matrix.  This is observed in \cite{Agrotis},
where one also finds the multi-Hamiltonian formulation:
\[ \overline{\pi}_{j}^\sharp \d h_{2-j}=\overline{\pi}_{j-1}^\sharp \d h_{3-j},   \quad k\in {\bf Z} \]
with $h_j=\frac{1}{j} \tr L^j$ for $j\not=0$ and $h_0=\ln(\det(L))$.

We would like to give now an intrinsic explanation for these
formulas, similar to the one given above for the Toda chain in
physical coordinates.

\subsubsection{Toda lattice in extended Flaschka coordinates}

Let us extend the Flaschka coordinates by considering a variable
$a_n$ defined by:
\[ a_n:={q_n}.\]
Then the transformation $(q_i,p_i)\mapsto (a_i,b_i)$ is a honest
change of coordinates. In these extended Flaschka coordinates, the
first Poisson tensor $\widetilde{\pi}_0$ is determined by:
\begin{align*}
\{a_i,b_i\}_0&=a_i,& &(i=1,\dots,n-1)&\\
\{a_i,b_{i+1}\}_0&=-a_i,& &(i=1,\dots,n-1)&\\
\{a_n,b_{n}\}_0&=1.&&&
\end{align*}
while the second Poisson $\widetilde{\pi}_1$ structure is given by:
\begin{align*}
\{a_i,a_{i+1}\}_1&=-a_ia_{i+1},& &(i=1,\dots,n-1)&\\
\{a_i,b_i\}_1&=a_ib_i,& &(i=1,\dots,n-1)&\\
\{a_n,b_{n}\}_1&=b_n \\
\{a_i,b_{i+1}\}_1&=-a_ib_{i+1},& &(i=1,\dots,n-1)&\\
\{b_i,b_{i+1}\}_1&=-a_i,& &(i=1,\dots,n-1).&
\end{align*}
In these coordinates, we still have the Nijenhuis tensor, relating
the various Poisson tensors in the hierarchy.

The submanifold $\R^{2n-1}\subset \R^{2n}$ defined by $a_n=0$ is a
Poisson submanifold for all Poisson tensors in the hierarchy, so
that the bi-Hamiltonian structure reduces to this submanifold, and
yields the bi-Hamiltonian formulation for the Toda lattice in
Flaschka coordinates. However, the tangent space to this submanifold
is not left invariant by the Nijenhuis operator $N$, and on
$\R^{2n-1}$ we do not have an induced PN-structure.

Another way of expressing these facts is to observe that the
involutive diffeomorphism $\phi:\R^{2n}\to\R^{2n}$ defined by:
\[ \phi(a_1,\dots,a_n,b_1,\dots,b_n)=(a_1,\dots,-a_n,b_1,\dots,b_n),\]
is a Poisson diffeomorphism for all Poisson structures. Hence, the
group $\Z_2=\{I,\phi\}$ acts by Poisson diffeomorphisms on
$\R^{2n}$, for all Poisson structures. It follows that its fix point
set, which is just $\R^{2n-1}$, has induced Poisson brackets (see
the Poisson Involution Theorem in \cite{FernandesVanhacke,Xu1}), and
these form the hierarchy in Flaschka coordinates.

\subsubsection{Toda lattice on a Lie algebroid}

We now consider the following bi-Ha\-mil\-to\-nian formulation on a
Lie algebroid.

We let $A=\R^{2n-1}\times\R^{2n}\to\R^{2n-1}$ be the trivial vector
bundle with fiber $\R^{2n}$. We denote by
$\{e_1,\dots,e_n,f_1,\dots,f_n\}$ a basis of global sections and we
let $(a_1,\dots,a_{n-1},b_1,\dots,b_n)$ be global coordinates on the
base. Now we define a Lie algebroid structure by declaring that the
bracket satisfies:
\[  [e_i,e_j]_A=[f_i,f_j]_A=[e_i,f_j]_A=0,\]
and that the anchor is given by:
\begin{align*}
\rho_A(e_i)&=\frac{\partial}{\partial a_i}\quad (i=1,\dots,n-1)&\rho_A(e_n)=0&\\
\rho_A(f_i)&=\frac{\partial}{\partial b_i}\quad(i=1,\dots,n).&
\end{align*}
Notice that $(A,[~,~]_A,\rho_A)$ is just the trivial extension of
Lie algebroids:
\[ \xymatrix{0\ar[r]&L_{\R}\ar[r]&A\ar[r]&T\R^{2n-1}\ar[r]&0}\]
where $L_{\R}$ denotes the trivial line bundle over $\R^{2n-1}$.

Now on $A$ we can define the following Poisson tensors:
\begin{align*}
    \pi_0&=\sum_{i=1}^{n-1} a_i e_i\wedge (f_i-f_{i+1})+e_n\wedge f_n\\
    \pi_1&=-\sum_{i=1}^{n-2} a_i a_{i+1}e_i\wedge e_{i+1}-a_{n-1}e_{n-1}\wedge e_n
      +\sum_{i=1}^{n-1} a_i e_i\wedge ( b_if_i-b_{i+1}f_{i+1})\\
 &\qquad\qquad+b_n e_n\wedge f_n-\sum_{i=1}^{n-1} a_i f_i\wedge f_{i+1}.
\end{align*}

These Poisson structures on $A$ cover ordinary Poisson structures on
the base $\R^{2n-1}$, which are just the Poisson structures
$\overline{\pi}_0$ and $\overline{\pi}_1$ of the Toda lattice, in
Flaschka coordinates.

Since $\pi_0$ is symplectic, the Poisson tensors on $A$ are
associated with a PN-algebroid structure. By our main theorem, they
give rise to an integrable hierarchy on $A$
\[ {\pi}_{j}^\sharp \d h_{2-j}={\pi}_{j-1}^\sharp \d h_{3-j},   \quad k\in {\bf Z} \]
with $h_j=\frac{1}{j} \tr N^j$ for $j\not=0$ and $h_0=\ln(\det(N))$,
covering an integrable hierarchy on the base
\[ \overline{\pi}_{j}^\sharp \d h_{2-j}=\overline{\pi}_{j-1}^\sharp \d h_{3-j}.   \quad k\in {\bf Z} \]
In this hierarchy the Hamiltonians differ by a factor of 2 relative to
the Hamiltonians in the multi-Hamiltonian formulation of the Toda lattice given by
Proposition \ref{propTodaMulti}. Although the hierarchy in the Lie
algebroid is generated by a Nijenhuis operator, it is well known
that this is not the case with the Toda lattice in the base
manifold.

\section{Acknowledgements}

This work was partially supported by POCI/MAT/58452/2004 and
CMUC/FCT.



\begin{thebibliography}{99}

\bibitem{Agrotis} M.~Agrotis and P.A.~Damianou, The modular hierarchy of the Toda lattice.
  \emph{Diff.~Geom Appl}, to appear.

\bibitem{BangouraKosmann} M.~Bangoura and Y.~Kosmann-Schwarzbach,
  Equation de Yang-Baxter dynamique classique et algebro{\"\i}des de Lie.
  \emph{C.~R.~Acad.~Sci.~Paris} \textbf{327} (1998), no. 8, 541--546.

\bibitem{CaWe} A.~Cannas da Silva and A.~Weinstein, \emph{Geometric
    Models for Noncommutative Algebras}, Berkeley Mathematics Lectures,
    vol.~\textbf{10}, American Math.~Soc.~, Providence, 1999.

\bibitem{EvLuWe} S.~Evens, J.-H.~Lu and A.~Weinstein, Transverse
  measures, the modular class and a cohomology pairing for Lie algebroids.
  \emph{Quart.~J.~Math.~Oxford (2)} \textbf{50} (1999), 417--436.

\bibitem{Fernandes1} R.L.~Fernandes, On the mastersymmetries and
  bi-Hamiltonian structure of the Toda lattice.
  \emph{J.~Phys.~A} \textbf{26} (1993), 3797--3803.

\bibitem{Loja} R.~L.~Fernandes, Lie algebroids, holonomy and
    characteristic classes. \emph{Adv.~in Math.~}\textbf{170} (2002),
    119--179.

\bibitem{LojaDamianou} R.~L.~Fernandes and P.~Damianou,
  Integrable hierarchies and the modular class. Preprint \emph{math.DG/0607784}.

\bibitem{FernandesVanhacke} R.L.~Fernandes and P.~Vanhaecke, Hyperelliptic
    Prym Varieties and Integrable Systems. \emph{Commun.~Math.~Phys.}~\textbf{221}
    (2001) 169--196.

\bibitem{Flaschka1}  H.~Flaschka, The Toda lattice. I. Existence of integrals.
  \emph{Phys.~Rev.~B (3)} \textbf{9} (1974), 1924--1925.

\bibitem{GraMarMich} J.~Grabowski, G.~Marmo and P.~Michor, Homology and
  modular classes of Lie algebroids.
  \emph{Ann.~Inst.~Fourier} \textbf{56} (2006), 69--83.

\bibitem{GraUrb} J.~Grabowski, P.~Urbanski, Lie algebroids and
Poisson-Nijenhuis structures.
  \emph{Reports~on~Mathematical~Physics} \textbf{40} (1997), no. 2, 195--208.

\bibitem{GraUrb2} J.~Grabowski, P.~Urbanski, Tangent and cotangent
lifts and graded Lie algebras associated with Lie algebroids.
\emph{Ann.~Global Analysis and Geometry} \textbf{15} (1997),
447--486.

\bibitem{Kosmann} Y.~Kosmann-Schwarzbach, The Lie bialgebroid of a
    Poisson-Nijenhuis manifold. \emph{Lett. Math. Phys.} \textbf{38} (1996), no. 4,
    421--428.

\bibitem{KosmannMagri0} Y.~Kosmann-Schwarzbach and F.~Magri, On the modular classes of Poisson-Nijenhuis manifolds. Preprint
\emph{math.SG/0611202}.

\bibitem{KosmannMagri} Y.~Kosmann-Schwarzbach and F.~Magri,
  Poisson-Nijenhuis structures.
  \emph{Ann.~Inst.~Henri Poincar\'e} \textbf{53} (1990), 35--81.

\bibitem{KosmannWeinstein} Y.~Kosmann-Schwarzbach and A.~Weinstein,
  Relative modular classes of Lie algebroids.
  \emph{C.~R.~Math.~Acad.~Sci.~Paris} \textbf{341} (2005), no. 8, 509--514.

\bibitem{Magri}  F.~{Magri},  P.~{Casati},   G.~{Falqui} and M.~{Pedroni}, Eight lectures on integrable
systems. In: Integrability of Nonlinear Systems (Y.
Kosmann-Schwarzbach et al. eds.), \emph{Lecture Notes in Physics}
\textbf{495} (2nd edition) (2004),  209–-250.



\bibitem{Weinstein} A.~Weinstein, The modular automorphism group of a
  Poisson manifold. \emph{J.~Geom.~Phys.~}\textbf{23} (1997), no.~3-4, 379--394.


\bibitem{Xu1} P.~Xu, Dirac submanifolds and Poisson involutions,
    \emph{Ann.~Sci.~\'Ecole Norm.~Sup.~}(4) \textbf{36} (2003), 403--430.

\end{thebibliography}
\end{document}